\documentclass[12pt]{amsart}
\usepackage{amssymb}
\usepackage{amsmath}
\usepackage[dvips]{graphicx}
\usepackage{psfrag}
\usepackage{hyperref}

\numberwithin{equation}{section}

\newcommand{\be}{\begin{equation}}
\newcommand{\ee}{\end{equation}}

\newcommand{\bone}{\bar{1}}
\newcommand{\btwo}{\bar{2}}

\newcommand\cS{{\mathcal S}}

\newtheorem{thm}{Theorem}
\newtheorem{cor}[thm]{Corollary}
\newtheorem{lem}[thm]{Lemma}
\newtheorem{defn}{Definition}
\newtheorem{rem}{Remark}

\makeatletter
\def\Ddots{\mathinner{\mkern1mu\raise\p@
\vbox{\kern7\p@\hbox{.}}\mkern2mu
\raise4\p@\hbox{.}\mkern2mu\raise7\p@\hbox{.}\mkern1mu}}
\makeatother

\title{Towards a human proof of Gessel's conjecture}
\author{Arvind Ayyer}
\address{Arvind Ayyer\\
Institut de Physique Th\'eorique\\
IPhT, CEA Saclay, and URA 2306, CNRS\\
91191 Gif-sur-Yvette Cedex, France}
\email{arvind.ayyer@cea.fr}
\date{\today}

\begin{document}

\begin{abstract}
We interpret walks in the first quadrant with steps
$\{(1,1),(1,0),(-1,0), (-1,-1)\}$ as a generalization
of Dyck words with two sets of letters. Using this language, we give a formal
expression for the number of walks in the steps above beginning and
ending at the origin. We give an explicit formula
for a restricted class of such words using a correspondance between
such words and Dyck paths.  
This explicit formula is exactly the same
as that for the
degree of the polynomial satisfied by the square of the area of cyclic
$n$-gons conjectured by Dave Robbins although the connection is a mystery. 
Finally we remark on another combinatorial problem in which the same
formula appears and argue for the existence of a bijection. 
\end{abstract}

\maketitle

\section{Introduction}
Ever since Gessel conjectured his formula for the number of walks in
the steps $\{(1,1),(1,0),(-1,0), (-1,-1)\}$ (which we will call Gessel
steps) starting and ending at the origin in $2n$ steps constrained to
lie in the first quadrant, there has been much interest in studying
lattice walks in the quarter plane. There have been conjectures for
lattice walks with Gessel steps terminating at other points \cite{pw},
as well as conjectures for the number of walks ending at the origin with
other sets of steps, most of which have been proven
\cite{bmm}. In a remarkable {\em tour de force}, Gessel's original  
conjecture has been finally proven using computer algebra techniques
\cite{kkz}. Even so, it is important to consider walks on the
quarter plane from a
human point of view because newer approaches tend to
open up  interesting mathematical avenues. 

In this article, we count a considerably restricted number of
walks with Gessel steps 
starting and ending at the origin by rephrasing the problem using
words with an alphabet consisting of four letters --- $1$, $2$, $\bone$ and
$\btwo$ which obey certain conditions. We first show that the
restatement of Gessel's conjecture in this context
 can be interpreted using Dyck paths. This gives a formal solution to
 the conjecture. Unfortunately, the solution is so formal as to be
 even computationally intractable!\footnote{Computing the $n$th term
   in the sequence involves $2n$ sums of binomial coefficients}.
We give a closed-form expression for the restricted problem and hope a
generalization of this method 
will give a better understanding of Gessel's conjecture. Admittedly
this result is a long way from a solution of the problem, but
one hopes that this technique can be generalized to obtain a complete proof
of Gessel's conjecture.

In Section \ref{sec:alph} we start with the preliminaries by defining
the alphabet and stating the main theorem. In Section \ref{sec:dyck},
we make the connection to Dyck paths and give a formal expression for
the number of walks beginning and ending at the origin using Gessel
steps. Section \ref{sec:pf} contains the proof which involves
summations of hypergeometric type. In principle, such sums can be
tackled by computer packages, but a certain amount of manipulation is
needed before they are summable. Lastly, we comment on related
problems in Section \ref{sec:open}.

\section{Gessel Alphabet} \label{sec:alph}
To rephrase the problem in the notation of formal languages, we need
some definitions. 

\begin{defn} \label{defn:galph}
The {\em Gessel alphabet} consists of a set of
letters $S = \{1,2,\dots\}$ with an order $<$ ($1<2<\dots$)
along with their complements which we denote
$\bar{S}=\{\bone,\btwo,\dots\}$. The order on the complement set is
irrelevant. 
\end{defn}

\begin{defn} \label{defn:gword}
Let $S = [n]$. Denote by $N_\alpha(w)$
the number of occurences of the letter $\alpha \in S \cup \bar S$ 
in the word $w$. (For example, $N_2(2 \btwo) = N_{\btwo}(2 \btwo) = 1,
N_1(2 \btwo) = 0$.) 
Then a {\em Gessel word} $w$ is a word such that every
prefix of the word satisfies 
$$\sum_{i=1}^k \left( N_{n+1-i}(w) -
N_{\overline{n+1-i}}(w) \right) \geq 0$$  
for each $k \in [1,n]$. 
\end{defn}

In words, this means that in each prefix, $n$ has to occur  more often
than $\bar{n}$, 
the number of occurrences of n and n-1 must be at least equal to 
 the number of occurrences of their barred
counterparts and so on. 
For example, $2 \bone$ is a
valid Gessel word but $1 \btwo$ is not. 

\begin{defn} \label{defn:cgword}
A {\em complete Gessel word} is a Gessel word $w$ where $N_i(w) =
N_{\bar i}(w)$ for all  letters $i \in S$. In other words, 
 the number of times the letter $i$ appears equals the number of
 times $\bar i$ appears for each $i \in [n]$.
\end{defn}

As an example, for $n=3$ both $3 \btwo 2 1 \bone \bar 3$ 
and $1 2 \bone \btwo 1$
are Gessel words but $2 1 \bar 3 2 3 \btwo \bone$ is not because the prefix
consisting of three letters fails the criterion in Definition
\ref{defn:gword}. Among the other two, the first one is a
complete Gessel word.

\begin{rem}
The number $G^{(d)}(n)$ of $d$ dimensional lattice walks in the first $2^d$-ant
with steps 
\be
\begin{split}
\{(1,\dots,1),(1,\dots,1,0),&\cdots,(1,0,\dots,0) \\
(-1,\dots,-1),(-1,\dots,-1,0),&\cdots,(-1,0,\dots,0)\} 
\end{split}
\ee
starting at the origin and returning
in $2n$ steps is the same as the number of complete Gessel words of length $2n$
in $d$ letters. None of these sequences seem to be present in
\cite{sloane} for dimensions higher than two and it would be
interesting to see if they are holonomic. Furthermore, none of these higher
dimensional sequences seem to have
the property of small factors which is present for the Gessel case.
\end{rem}

$G^{(2)}(n)$ is conjectured by Gessel to be given by
the closed form expression
\be \label{gconj}
16^n \frac{(5/6)_n (1/2)_n}{ (2)_n (5/3)_n},
\ee
where
\be
(a)_n = a (a+1) \dots (a+n-1)
\ee
is the Pochhammer symbol or  rising factorial. The first few terms
are the sequence A135404 
in \cite{sloane}. For the remainder of the paper, we implicitly assume
$d=2$ and omit the superscripts in defining various constrained Gessel
numbers. 

We express the number of walks of length $2n$
 as a number triangle based on the number of times $2$ and
 $\btwo$ appear increasing from left to right.
\be
\begin{array}{ccccccccc}
 & & & &1 & & & & \\ 
 & & &1 & &1 & & & \\ 
 & &2 & &7 & &2 & & \\ 
 &5 & &37 & &38 & &5 & \\ 
14 & &177 & &390 & &187 & &14 
\end{array}
\ee

One immediately notices that the leftmost and rightmost entries are
the Catalan numbers. This is because the number of complete Gessel
words with $N_2(w) = 0$ ($N_1(w)= 0$) in a word $w$ of length
$2n$ is
 in immediate bijection with the number of Dyck paths ending at
 $(2n,0)$ because $N_1(x) > N_{\bone}(x)$ ($N_2(x) > N_{\btwo}(x)$) for
 each prefix $x$ of the word $w$.

What is more interesting, and the main result of the paper is the
next-to-rightmost sequence beginning $1,7,38,187$. Strangely
enough, this sequence is already present in the OEIS as A000531
\cite{sloane}. It turns out to be exactly the one conjectured by Dave Robbins
\cite{robb1} to be the degree of the polynomial satisfied by $16K^2$,
where $K$ is the area of a cyclic $n$-gon and proved in \cite{fp,mrr}. As
far as we know, this 
result is a coincidence without any satisfactory explanation.
For a recent review of the
subject, see \cite{pak}. This is also related to Simon Norton's
conjecture on the same page in the OEIS. We comment on this in Section
\ref{sec:open}. 

\begin{thm} \label{thm:gws}
The number of complete Gessel words $G_1(n)$ in two letters with $n-1$ $2$'s and
$\btwo$'s, and one $1$ and $\bone$, is given by
\be \label{cn}
G_1(n) = \frac{(2n+1)}{2}\binom{2n}{n} - 2^{2n-1}.
\ee
\end{thm}

The proof uses the idea that the number of Gessel words with $n_2$
$2$'s and $\btwo$'s and $n_1$ $1$ and $\bone$'s can be
calculated using a bijection with Dyck paths. The answer can be
written as a sum of 
products of expressions counting the number of Dyck paths between two
different heights. The summation can be done explicitly
 when $n_1=1$.

\section{Complete Gessel words and Dyck paths} \label{sec:dyck}
We consider Dyck paths to be paths using steps
$\{(1,1),(1,-1)\}$ starting at the origin, staying on or above the
$x$-axis and ending on the $x$-axis. 
In this section we exhibit a bijection between complete Gessel words
(the counting of which is stated by the conjecture of Gessel)
and a set of 
restricted Dyck paths which will be useful in the proof of
Theorem~\ref{thm:gws}.

\begin{defn}
Let $P=(P_1,\dots,P_{m})$ be an increasing list of positive integers
and $H=(H_1,\dots,H_{m})$ be a list of nonnegative integers of the
same length. We define a {\em $(P,H)$-Dyck path} to be a Dyck path of
length greater than $m$ which satisfies the constraint that between
positions $P_i$ and $P_{i+1}$ (both inclusive), the ordinate of the
path is greater than or equal to $H_i$ for $i=1,\dots,m-1$.
\end{defn}

Notice that this forces the ordinates of the path at positions $P_i$
to be greater than or equal to the heights $\max\{H_{i-1},H_i\}$.

We now associate to every complete Gessel word $w$ in two letters lists
$P$ and $H$ using the following algorithm. 
\begin{enumerate}
\item Construct the list $S$ of length $2n_1$ of letters $1$ or $\bone$ as
  they occur in the word. 
\item From the list $S$, construct the list $T$ by replacing $1$ by
  $1$ and $\bone$ by $-1$.
\item Construct the list $\tilde P$ whose elements are positions of
  the letter $S_i$ in $w$. Similarly, construct the elements of the
  list $P$ as $\tilde P_i -i$.
\item Finally, each element of the list $H$ is given by
\be \label{hts}
H_i = \max\Big\{-\sum_{k=1}^i T_k,0\Big\}.
\ee
\end{enumerate}
Clearly, $S$ and $H$ determine each other and similarly, so do $P$ and
$\tilde P$. Therefore one can also associate a complete Gessel word to
a $(P,H)$-Dyck path and vice versa.
As an example, consider the Gessel word $w=2 \bone 2 1 \btwo
\btwo$. For this word, $S=(\bone,1), T=(-1,1), \tilde P = (2,4), P=(1,2)$ and
$H=(1,0)$. Also, given this $P$ and $H$, there is exactly one such $(P,H)$-Dyck
path of length four, namely
$(\nearrow,\nearrow,\searrow,\searrow)$, just as $w$ is the only 
complete Gessel word of length six with $\bone$ at position two at $1$
at position four.

\begin{lem} \label{lem:bij}
Complete Gessel words of length $2(n_1+n_2)$ in two letters with
positions of $1$ and $\bone$ given by the lists $\tilde P,S$ are in bijection
with $(P,H)$-Dyck paths of length $2n_2$ where the pairs of lists
$(\tilde P,S)$ and $(P,H)$ are
related by the algorithm described above.
\end{lem}

\begin{proof}
Starting with the complete Gessel word, one replaces each occurence of
the letter $2$
by the step $(1,1)$ and that of  $\btwo$ by the step $(1,-1)$. The
constraint defining the 
$(P,H)$-Dyck path is simply another way of expressing the inequality in
Definition~\ref{defn:gword}. 
\end{proof}

One could generalize this bijection to include paths not ending on the
$x$-axis and Gessel words which are not complete, but this is
sufficient for our purposes.

One of the main tools in the proof of Theorem~\ref{thm:gws} is an 
expression for the number of Dyck paths between two different heights,
which can be readily obtained from the reflection principle \cite{co}.

\begin{lem} \label{lem:dyck}
The number of Dyck paths $a_{i,j}(k)$
 that stay above the $x$-axis starting at the
position $(0,i)$ and end at position $(k,j)$ is given by
\be \label{aij}
a_{i,j}(k) = \begin{cases}
\displaystyle \binom{k}{(k+i-j)/2} -\binom{k}{(k+i+j)/2+1} & \\
 & \hspace{-2cm} \text{if } (k+i+j)  \equiv 0 \mod 2, \\
0 & \hspace{-2cm} \text{if } (k+i+j) \equiv 1 \mod 2.
\end{cases}
\ee
\end{lem}

We now use the bijection in Lemma~\ref{lem:bij} and the formula in
Lemma~\ref{lem:dyck}
to write an expression for the number of complete Gessel
words of length $2n$ for fixed positions of $1,\bone$.

\begin{lem} \label{lem:gwn2}
Let us fix the positions of $n_1$ $1,\bone$ by the lists $S,\tilde P$.
Calculate the lists $T$ and $H$ by the algorithm above and 
let $G_{n_1}(S,\tilde P;2n)$ denote the number of such complete Gessel
words. Then
\be \label{cpt}
\begin{split}
&G_{n_1}(S, \tilde P;2n) = 
\displaystyle \sum_{k_1=\delta_{(1-T_1)/2,0}}^{\tilde P_1-1}
  \sum_{k_2=H_2}^{\tilde P_2-1} \dots  \sum_{k_i=H_i}^{\tilde P_i-1} \dots \sum_{k_{2n_1}=H_{2n_1}}^{\tilde P_{2n_1}-1} \\
&\displaystyle a_{0,k_1}(\tilde P_1-1) \; a_{k_{2n_1},0}(2n-\tilde P_{2n_1}) \prod_{i=2}^{2n_1-1}
a_{k_{i-1}-H_i,k_i-H_i}(\tilde P_i-\tilde P_{i-1}-1),
\end{split}
\ee
where the lower index of the sum $k_1$ depends on the first element of
the list $T$. 
\end{lem}

\begin{proof}
The proof is straightforward, using the bijection of
Lemma~\ref{lem:bij} to rewrite each
Gessel word with the positions of $1,\bone$ given by the lists
$S, \tilde P$ as a Dyck path with heights at the points $P_i$ (given
by $k_i$) being not less than $H_i$
and then the reflection principle in Lemma~\ref{lem:dyck} to count the
number of paths between position $P_{i-1}$ and $P_i$ for each $i$.
\end{proof}

\begin{cor} \label{cor:cat}
For a given configuration of $1,\bone$, replace each $+1$ in $T$ by
an upward Dyck step and each $-1$ by a downward Dyck step. If the
whole of $T$ forms an legal Dyck path, then $G_{n_1}(S, \tilde P;2n) =
C_{n_1}$, the $n_1$th Catalan number independent of the list $P$.
\end{cor}

\begin{proof}
Whenever the above condition is satisfied, $H_i=0$ for all $i$, which
means we simply count the number of Dyck paths of length $2n_1$ in
\eqref{cpt} by definition.
\end{proof}

Now we obtain a formula for the number of complete Gessel words with
$n_1$ $1,\bone$'s using
Lemma~\ref{lem:gwn2} and writing down all possibilities for $\tilde P$ and $S$.
The number of ways of writing all possible $\tilde P$'s is simply
$\binom{2n}{2n_1}$ because one has to choose $2n_1$ positions out of
$2n$ positions. For each $\tilde P$, one has to choose $n_1$ positions for
$1$ and $\bone$ each and therefore the number of such ways is
$\binom{2n_1}{n_1}$.

Let us form the set 
\be
\cS = \left\{(S,\tilde P) \Biggl\vert \begin{array}{l}
S  \text{ is an ordered list of $n_1$ $1$'s and  $n_1$ $\bone$'s.} \\
\tilde P \text{ is an increasing list of} \\
 \; \text{ $2n_1$ positions between $1$ and $2n$,}
\end{array} \right\},
\ee
whose cardinality is
\be
\binom{2n}{2n_1} \binom{2n_1}{n_1} = \frac{(2n)!}{(n_1)!^2 (2n-2n_1)!}.
\ee
Therefore the number of complete Gessel words with exactly $n_1$
$1,\bone$'s is given by
\be \label{cgess}
G_{n_1}(n) = \sum_{(S,\tilde P) \in \cS} G_{n_1}(S, \tilde P;2n),
\ee
and the number of complete Gessel words in $2n$ letters is
\be
G(n) = \sum_{n_1=0}^n G_{n_1}(n).
\ee
Showing that $G(n)$ is equal to the expression \eqref{gconj} would be
the ultimate (and possibly hopeless) aim of this line of approach.

We now have all the ingredients necessary to prove
Theorem~\ref{thm:gws} which corresponds to  the special case
$n_1=1$. Before we go on to 
the proof, however, we make some observations about complete Gessel
words with exactly one $1$ and $\bone$.
Let $d_{i,j}$ be the number of times there is an $1$ or a $\bone$ at
position $i$ and its counterpart at position $j$. Then we draw the
following triangle for a specific $n$,
\be
\begin{array}{ccccccccc}
 & & & &d_{1,2n} & & & & \\ 
 & & &d_{1,2n-1} & &d_{2,2n} & & & \\ 
 & &d_{1,2n-2} & &d_{2,2n-1} & &d_{3,2n} & & \\ 
 &\Ddots& & & & & &\ddots & \\ 
d_{1,2} & &d_{2,3} & &\cdots & &d_{2n-2,2n-1} & &d_{2n-1,2n}.
\end{array}
\ee
For $n=3$, the triangle is
\be
\begin{array}{ccccccccc}
 & & & &2 & & & & \\ 
 & & &2 & &2 & & & \\ 
 & &2 & &3 &  &2 & & \\ 
 &2 & &3 & &3 & &2 & \\ 
2 & &4 & &3 & &4 & &2, 
\end{array}
\ee
and for $n=4$, the triangle is
\be
\begin{array}{ccccccccccccc}
 & & & & & &5 & & & & & &\\ 
 & & & & &5 & &5 & & & & &\\ 
 & & & &5 & &7 & &5 & & & &\\ 
 & & &5 & &7 & &7 & &5 & & &\\ 
 & &5 & &8 & &7 & &8 & &5 & & \\
 &5 & &8 & &8 & &8 & &8 & &5 &\\ 
5 & &10 & &8 & &10 & &8 & &10 & &5. 
\end{array}
\ee
It is clear that if one stacks these triangles on top of one another,
one gets a generalization of Pascal's pyramid, where each layer $n$
fits in the vacancies of the layer $n-1$ above it.

\begin{rem}
We note some properties of these triangles.
\begin{enumerate}
\item The sum of the entries in the triangle are precisely what we
  claim are given by \eqref{cn}.

\item One notices immediately that the extremal columns are Catalan numbers
$C_{n-1}$. This follows immediately from Corollary~\ref{cor:cat} and
  the fact that a Gessel word cannot begin with $\bone$ or $\btwo$, or
  end with a $1$ or $2$. The even entries in the last row are $2
  C_{n-1}$. 

\item Every number in the interior of the triangle occurs $4k$ times
  for $k$ a positive integer. Furthermore, they are organized as
  rhombus-shaped blocks of size four.
This turns out to be true for all
$n$. We will need this fact in the proof later and we state it as
Lemma~\ref{lem:diamond}. 
\end{enumerate}
\end{rem}

\section{Proof of Theorem~\ref{thm:gws}} \label{sec:pf}

One simply has to analyze all possibilities of occurences of $1$ and
$\bone$ case by case. Suppose $1$ occurs at
position $i$ and $\bone$ occurs at position $j$ in a word of length
$2n$ and $i<j$. Then by Corollary~\ref{cor:cat}, the number of such Gessel words
is $C_{n-1}$. The number of possibilities
of $i,j$ such that $1 \leq i<j \leq 2n$ is
$n(2n-1)$. Therefore, the number of Gessel words where the $1$
occurs before the $\bone$ is 
\be \label{gw1bone}
(2n-1)\binom{2n-2}{n-1}.
\ee

We now use Lemma~\ref{lem:gwn2} to count the number of words where
$\bone$ occurs at site $i$ before $1$ at site $j$,
\be 
\begin{split}
G_1([i,j],&[-1,1];2n) \\= &\sum_{k_1=1}^{i-1} \sum_{k_2=0}^{j-1}
\, a_{0,k_1}(i-1)
 \, a_{k_1-1,k_2-1}(j-i-1) \, a_{k_2,0}(2n-j),
\end{split}
\ee
which, using \eqref{aij} gives
\be \label{c1}
\begin{split}
G_1([i,j],&[-1,1];2n) =\sum_{k_1=1}^{i-1} \sum_{k_2=0}^{j-1}\; C^{(i-1+k_1)/2}_{(i-1-k_1)/2} \;
C^{(2n-j+k_2)/2}_{(2n-j-k_2)/2} \\
&\left[
\binom{j-i-1}{(j-i-1+k_1-k_2)/2} - \binom{j-i-1}{(j-i-1+k_1+k_2)/2} \right]
\end{split}
\ee
where $C^m_n$ is the Catalan triangle number given by $\frac{(m-n+1)}{(m+1)}
\binom{m+n}{n}$ for $0 \leq n \leq m, m\geq 0$.

We now use the following result to simplify calculations. The proof of
this assertion is easily verified by expanding \eqref{c1} and noting
that the answer is the same when $i$ is replaced by either $2i$ or
$2i+1$ and similarly for $j$.
\begin{lem} \label{lem:diamond}
For $1 \leq i < j \leq n-1$,
\be
\begin{split}
&G_1([2i,2j],[-1,1];2n) = G_1([2i,2j+1],[-1,1];2n) \\
&= G_1([2i+1,2j],[-1,1];2n) =
G_1([2i+1,2j+1],[-1,1];2n).
\end{split}
\ee
\end{lem}

Then the total number of Gessel words with an $\bone$ preceding an
$1$ is given by
\be \label{sumc1}
\begin{split}
\sum_{i=1}^{2n-1} \sum_{j=i+1}^{2n} G_1([i,j],[-1,1];2n) = &4
\sum_{i=1}^{n-2} \sum_{j=i+1}^{n-1} G_1([2i,2j],[-1,1];2n) \\
&+ \sum_{i=1}^{n-1} G_1([2i,2i+1],[-1,1];2n) \\
&= 4(S_2 - S_3) + S_1.
\end{split}
\ee
where we have split the sum in three parts, with 
\be \label{s1}
S_1 = \sum_{i=1}^{n-1} G_1([2i,2i+1],[-1,1];2n).
\ee
The remainder in \eqref{sumc1} we split using \eqref{c1}, and using the
variables $r=(2i-k_1-1)/2,s=(2n-2j-k_2)/2$, as
\be
\sum_{i=1}^{n-2} \sum_{j=i+1}^{n-1} G_1([2i,2j],[-1,1];2n) = S_2 - S_3
\ee
where
\be \label{s2s3}
\begin{split}
S_2 &= \sum_{i=1}^{n-2} \sum_{j=i+1}^{n-1} \sum_{r=0}^{i-1}
\sum_{s=0}^{n-j-1} C^{2i-r-1}_r C^{2n-2j-s}_s
\binom{2j-2i-1}{2j+s-n-r-1}, \\
S_3 &= \sum_{i=1}^{n-2} \sum_{j=i+1}^{n-1} \sum_{r=0}^{i-1}
\sum_{s=0}^{n-j-1} C^{2i-r-1}_r C^{2n-2j-s}_s
\binom{2j-2i-1}{n-s-r-1}. 
\end{split}
\ee
We now estimate these three sums in turn.

\subsection{The sum $S_2$}
Replacing $r \to i-1-r$ and $s \to n-j-s-1$,
substituting $k=j-i$ and rearranging the variables, we get 
\be 
S_2 = \sum_{r=0}^{n-3} \sum_{k=1}^{n-r-2} \sum_{i=r+1}^{n-1-k}
\sum_{s=0}^{n-i-k-1} C^{i+r}_{i-r-1} C^{n-k-i+s+1}_{n-k-i-s-1}
\binom{2k-1}{k-s+r-1}.
\ee
Now replace $k \to k-1, i \to i-r-1$ to get
\be 
S_2 = \sum_{r=0}^{n-3} \sum_{k=0}^{n-r-3} \sum_{i=0}^{n-r-k-3}
\sum_{s=0}^{n-r-k-i-3} C^{i+2r+1}_{i} C^{n-k-i-r+s-1}_{n-k-i-r-s-3}
\binom{2k+1}{k-s+r}.
\ee
We now replace the $r$ variable by $u=k+r$. Notice that the binomial
coefficient term is independent of $i$ for which we use the identity
\be \label{catid}
\sum_{i=0}^{C} C^{i+A}_{i} C^{B-i}_{C-i} =
C^{A+B+1}_{C},
\ee
which means we are left with
\be
\begin{split}
S_2 &= \sum_{u=0}^{n-3} \sum_{k=0}^{u} \sum_{s=0}^{n-u-3}
C^{u+n+s-2k+1}_{n-u-s-3} \binom{2k+1}{u-s}\\
&= \sum_{k=0}^{n-3} \sum_{u=0}^{n-k-3} \sum_{s=0}^{n-u-k-3}
C^{u+n+s-k+1}_{n-u-s-k-3} \binom{2k+1}{u+1-s}.
\end{split}
\ee
Let $A=n-k-3, v=u-s$ and $v'=s-u$. Then one easily verifies that
\be \label{splitsum}
\sum_{u=0}^A \sum_{s=0}^{A-u} = \sum_{v=0}^A \sum_{\substack{u=v \\
(s=u-v)}}^{A/2+ v/2} +
\sum_{v'=0}^A \sum_{\substack{s=v' \\ (u=s-v')}}^{A/2+v'/2} -
\sum_{\substack{s=0 \\ (u=s)}}^{A/2}
\ee 
The binomial coefficient is independent of
$u$ in the first sum and of $s$ in the remaining two and hence the
innermost sum can be done using the identity
\be
\sum_{s=v}^{n/2} C^{B+2s}_{n-2s} = \binom{B+n-1}{n-2v}.
\ee

This reduces the sum (after a change of variables) to
\be
\begin{split}
S_2 = \sum_{k=0}^{n-3} \sum_{v=0}^{k} &\binom{2k+3}{k-v}
\binom{2n-2k-4}{n-k-v-2} \\ &- \sum_{k=0}^{n-3} \binom{2k+1}{k+1}
\binom{2n-2k-3}{n-k-3} 
\end{split}
\ee
These sums are handled as special cases of the Chu-Vandermonde
identity to yield
\be \label{s2ans}
S_2 = \frac{n+2}{4} \binom{2n}{n} - 3 \cdot 2^{2n-3},
\ee
which appears as sequence A045720 \cite{sloane} because it is the
threefold convolution of the sequence $a_n=\binom{2n+1}{n+1}$.

\subsection{The sum $S_3$}
\be
S_3 = \sum_{i=1}^{n-2} \sum_{j=i+1}^{n-1} \sum_{r=0}^{i-1}
\sum_{s=0}^{n-j-1} C^{2i-r-1}_r C^{2n-2j-s}_s
\binom{2j-2i-1}{n-s-r-1},
\ee
which after replacing $i \to i-r-1$ and subsequently $j \to j-r-i-2$
and rearranging becomes
\be
S_3 = \sum_{j=0}^{n-3} \sum_{i=0}^{n-j-3} \sum_{s=0}^{n-j-i-3}
\sum_{r=0}^{n-j-i-s-3} C^{2i+r+1}_r C^{n-j-r-i+s-1}_{n-j-i-r-s-3}
\binom{2j+1}{s+j+i+2}.
\ee
We now use \eqref{catid} to do the $r$ sum and get
\be
S_3 = \sum_{j=0}^{n-3} \sum_{i=0}^{n-j-3} \sum_{s=0}^{n-j-i-3}
C^{n-j+i+s+1}_{n-j-i-s-3} \binom{2j+1}{s+j+i+2}.
\ee
Now, replacing $s$ by $k=i+s$, we get 
\be
\begin{split}
S_3 &= \sum_{j=0}^{n-3} \sum_{k=0}^{n-j-3} \sum_{i=0}^{k}
C^{n-j+k+1}_{n-j-k-3} \binom{2j+1}{k+j+2} \\
&= \sum_{j=0}^{n-3} \sum_{k=0}^{n-j-3} 
(k+1) C^{n-j+k+1}_{n-j-k-3} \binom{2j+1}{k+j+2}\\
&= \sum_{k=0}^{n-3} \sum_{j=0}^{n-k-3} 
(k+1) C^{n-j+k+1}_{n-j-k-3} \binom{2j+1}{k+j+2}.
\end{split}
\ee
We now use the identity
\be \label{catid2}
\sum_{j=C}^B C^{A-j}_{B-j} \binom{2j+1}{j-C} = \binom{A+B+2}{B-C},
\ee
for the $j$ sum to get 
\be \label{s3ans}
\begin{split}
S_3 &= \sum_{k=0}^{n-3}(k+1) \binom{2n}{n-4-2k},\\
 &= \sum_{k=0}^{(n-4)/2}(k+1) \binom{2n}{n-4-2k},\\
 &= \frac12 \sum_{k=0}^{(n-4)/2}(2k+4)
 \binom{2n}{n-4-2k}-\sum_{k=0}^{(n-4)/2} \binom{2n}{n-4-2k}\\
 &= \frac n2
 \binom{2n-2}{n-4}-\sum_{k=0}^{(n-4)/2} \binom{2n}{n-4-2k}\\
 &= \frac n2
 \binom{2n-2}{n-4}-2^{2n-2} + \frac{(2n)!  (3n^2+n+2)}{2 n! (n+2)!}.
\end{split}
\ee

\subsection{The sum $S_1$}
\be
S_1 = \sum_{i=1}^{n-1} \sum_{r=1}^i \sum_{s=1}^i C^{i-1+r}_{i-r}
C^{n-i+s-1}_{n-i-s} \left[ \binom{0}{r-s} - \binom{0}{r+s-1} \right].
\ee
The first term forces $r=s$ and the second term is identically zero
because $r+s \geq 2$. This means we are left with
\be \label{s1ans}
\begin{split}
S_1 &= \sum_{i=1}^{n-1} \sum_{r=1}^i  C^{i-1+r}_{i-r}
C^{n-i+r-1}_{n-i-r} \\
&=  (n-1) C_{n-1}.
\end{split}
\ee

Thus the total number of Gessel words where an $\bone$ occurs before
an $1$ defined in \eqref{sumc1} is given, using
\eqref{s2ans},\eqref{s3ans} and \eqref{s1ans}, by
\be \label{gwbone1}
4(S_2-S_3) + S_1 = \frac{(n^3+4n^2+5n+2) (2n)!}{2 n! (n+2)!} - 2^{(2n-1)},
\ee
and therefore, the total number of complete Gessel words is
\be
G_1(n) = 4(S_2-S_3) + S_1 + (2n-1) \binom{2n-2}{n-1} =
\frac{(2n+1)}{2}\binom{2n}{n} - 2^{2n-1},
\ee
which is exactly the same expression as \eqref{cn}. 
\hfill $\blacksquare$

\section{Remarks} \label{sec:open}
This section is intended to be speculative in nature and
consequently, the statements are unproven as far as we know, though
not necessarily very 
deep. 
In 2001, Simon Norton made the following conjecture in A000531 \cite{sloane}.
\begin{center}
{\em A conjectured definition: Let $0 < a_1 < a_2 <...<a_{2n} < 1$. Then how
many ways are there in which one can add or subtract all the $a_i$ to
get an odd number. For example, take $n = 2$. Then the options are
$a_1+a_2+a_3+a_4 = 1$ or $3$; one can change ths sign of any of the $a_i$'s
and get 1; or $-a_1-a_2+a_3+a_4 = 1$. That's a total of 7, which is the
2nd number of this sequence.}
\end{center}

We want to connect this conjecture to Theorem \ref{thm:gws}. 
Before that, we need some preliminaries. 
One can
represent every equation of the form $\pm a_1 \dots \pm a_{2n}=1$ as a
$2n$-tuple of $+,-$ symbols.
Let us replace every $-$ by a 0 and every $+$ by a 1. Then,
one can represent all possible ways
of ordering the $+$'s and $-$'s by binary words of length $2n$.

Let $w$ be such a binary word. Then define $n_1(w)$ to be number of
1's in w. 
Also define $n_{10}(w)$ to be the number of occurences of distinct 10
subwords in $w$. For example, $n_{10}(1110) = 1$ and $n_{10}(0110000)=2$.
We now form the multiset $S$, where each word $w$ occurs 
\be
m(w) = \bigg\lfloor \frac{n_1(w) -n_{10}(w)}2 \bigg\rfloor
\ee 
times. Note that if
$m(w)$ is zero or negative, it never appears. Then, it seems that the
cardinality of $S$ is the same as the conjecture in the sequence. 
Moreover there is a
bijection from the $\pm$ notation to the binary notation. 
This means that the number of times a binary word appears in $S$ seems
to be the same as the number of positive odd integers in the right
hand side of the equation corresponding to the same binary word which
admit solutions. We give a concrete example in Table~\ref{tab:n2}.
\begin{table}[h]
\begin{tabular}{|c|c|c|c|c|c|} 
\hline
$\pm$ &  Odd integer &     Binary  &  $n_1(w)$ &
  $n_{10}(w)$ & $m(w)$ \\
 word &   sums & word & & & \\
\hline                           
$++++$ &           1,3 &            1111 &   4 &      0   &            2\\
$+++-$ &           1   &            1110 &   3 &      1   &            1\\
$++-+$ &           1   &            1101 &   3 &      1   &            1\\
$+-++$ &           1   &            1011 &   3 &      1   &            1\\
$-+++$ &           1   &            0111 &   3 &      0   &            1\\
$--++$ &           1   &            0011 &   2 &      0   & 1\\
\hline
\end{tabular}
\caption{All allowed possibilities for $n=2$.  \label{tab:n2}}
\end{table}

The connection between the two problems is as follows. For each fixed
number $n_1$ 
of $+$ signs from 2 to $2n$, count only those sums in which all
possible $\binom{2n}{n_1}$ combinations give rise to that sum and add
them up. This number is precisely the same as the number of Gessel
words stated in Theorem~\ref{thm:gws} in which the 1 precedes the
$\bar{1}$.  The formula for the number of such Gessel words is given
by \eqref{gw1bone}. 
If one considers the set of only those $\pm$ words for fixed $n_1$
such that a number strictly smaller than $\binom{2n}{n_1}$ contribute,
then this set is equinumerous with the
Gessel words stated above in which the $\bar{1}$ precedes the 1 and is
given by \eqref{gwbone1}. This
leads us to conjecture the presence of a bijection between the
multiset $S$ and the number of complete Gessel words with exactly one
$1$ and $\bone$.

\begin{table}[h]
\begin{tabular}{|c|c|c|c|c|}
\hline
Number of $+$ and & Sum=1 & Sum=3 & Sum=5 & Sum=7 \\ 
and $-$ signs & & & & \\
\hline
$8+$ & 1  & 1  & 1  & 1 \\
$7+,1-$ & 8  & 8  & 8  &  \\
$6+,2-$ & 28  & 28  & 1  &  \\
$5+,3-$ & 56  & 8  &   &  \\
$4+,4-$ & 28  & 1  &   &  \\
$3+,5-$ & 8  &   &  &  \\
$2+,6-$ & 1  &   &   &  \\
\hline
\end{tabular}
\caption{The number of words for a fixed number of $+$ and $-$ signs
  and fixed sum in the case $n=4$. \label{tab:n4}}
\end{table}

For example, there are 6 complete Gessel words for $n=2$ where the 1
precedes the $\bone$. From Table~\ref{tab:n2}, one sees that all
possible terms contribute when we have either $4+$ or $3+,1-$
signs. There are two possibilities for the former (when the sums are 1 and 3)
and four for the latter (when the sum is 1). Similarly, there is only
one complete Gessel word for $n=2$ where $\bone$ precedes 1, which is
given by $2 \bone 2 1 \btwo \btwo$ and for $2+,2-$ signs, there
are 6 possible words, but only one contributes.


For any fixed $n_1$
and any fixed odd integer sum, the number of words which allow this 
seem to be of the form $\binom{2n}{k}$ where k varies from $0$ to
$n-1$! We illustrate this via another concrete example in
Table~\ref{tab:n4}. Notice that the only integers appearing in the
table are the binomial coefficients  $\binom{8}{k}$ with $k=0,1,2$ or 3.
Another observation is that if one draws lines of $45^\circ$ starting
from the first column in Table~\ref{tab:n4} and looks at the diagonal
columns, one finds the pattern,
\be \label{n4col}
\begin{array}{|c|c|c|c|c|c|c|}
1 & 1 & 1 & 1 & 1 & 1 & 1 \\
  & 8 & 8 & 8 & 8 & 8  &\\
  &   & 28 & 28 & 28 & & \\
  &   &    & 56 &   &   & 
\end{array}
\ee
from which it is clear that each of these diagonal columns in
\eqref{n4col} starts with
$\binom{2n}{0}$ with 
subsequent values of the lower index increasing by 1. The 
first four columns above correspond exactly to the Gessel words where
$1$ precedes $\bone$ is the sum of the entries is precisely $(2n-1)
\binom{2n-2}{n-1}$ with $n=4$. This pattern persists up until $n=6$.

\section*{Acknowledgements}
The author would like to thank Doron Zeilberger and Tewodros
Amdeberhan for useful discussions, Lara Pudwell and Eric Rowland for
many helpful comments on an earlier version of the draft, and Simon
Norton for discussing his conjecture.

\end{document}